\documentclass[12pt,reqno]{amsart}

\usepackage{xcolor}
\newtheorem{theorem}{Theorem}
\newtheorem{lemma}{Lemma}
\newtheorem{proposition}{Proposition}
\newtheorem{remark}{Remark}

\newtheorem{corollary}{Corollary}

\newcommand{\beq}{\begin{equation}}
\newcommand{\eeq}{\end{equation}}
\newcommand{\beqna}{\begin{eqnarray*}}
\newcommand{\eeqna}{\end{eqnarray*}}
\newcommand{\beqn}{\begin{equation*}}
\newcommand{\eeqn}{\end{equation*}}
\newcommand{\bp}{\begin{proof}}
\newcommand{\ep}{\end{proof}}
\newcommand{\bprop}{\begin{proposition}}
\newcommand{\eprop}{\end{proposition}}
\newcommand{\bt}{\begin{theorem}}
\newcommand{\et}{\end{theorem}}
\newcommand{\bex}{\begin{Example}}
\newcommand{\eex}{\end{Example}}
\newcommand{\bc}{\begin{corollary}}
\newcommand{\ec}{\end{corollary}}
\newcommand{\bl}{\begin{lemma}}
\newcommand{\el}{\end{lemma}}
\newcommand{\br}{\begin{remark}}
\newcommand{\er}{\end{remark}}
\newcommand{\eps}{\varepsilon}
\newcommand{\R}{{\mathbb R}}
\newcommand{\RR}{{\mathcal R}}

\begin{document}

\title
[Distributions of Sections and Projections]
{Distribution functions of sections and projections of convex bodies}

\author{Jaegil Kim, Vladyslav Yaskin, and Artem Zvavitch}\thanks{The first and second named authors are supported in part by NSERC. The third named author is supported in part by  the U.S. National Science Foundation Grant DMS-1101636 and by the  Simons Foundation. Part of this work was done during the second and third authors' stay at the Mathematisches Forschungsinstitut Oberwolfach. The authors gratefully acknowledge the Institute's hospitality.}

\address{Department of Mathematical and Statistical Sciences, University of Alberta, Edmonton, Alberta T6G 2G1, Canada} \email{jaegil@ualberta.ca}

\address{Department of Mathematical and Statistical Sciences, University of Alberta, Edmonton, Alberta T6G 2G1, Canada} \email{vladyaskin@math.ualberta.ca}

\address{Department of Mathematics, Kent State University,
Kent, OH 44242, USA} \email{zvavitch@math.kent.edu}

\date{\today}


\begin{abstract}
Typically, when  we are given the section (or projection) function of a convex body, it means that in each direction we know the size of the central section (or projection) perpendicular to this direction.
Suppose now that we can only get the information about the sizes of sections (or projections), and not about the corresponding directions.
In this paper we study to what extent the distribution function of the
areas of central sections (or projections) of a convex body can be used to derive some information about the body, its volume, etc.
\end{abstract}

\maketitle

\section{Introduction}

Let $K$ be a convex body in $\R^n$. What can we say about the body  $K$, if we know the areas of all its central sections or projections in every direction?  Such questions are typical in Geometric Tomography. Particular examples include questions on the unique determination, volume comparison problems, etc. In this paper we study  similar questions under much weaker assumptions. Instead of the pointwise knowledge of the areas of sections or projections, let us assume that we merely have their distribution functions.

Let $\sigma$ be the Haar probability measure on $S^{n-1}$.
If $K \subset \R^n$ is a convex body that contains the origin in its interior, define the functions $S_K(t)$ and $\Pi_K(t): \R^+ \to [0,1]$ by
$$
S_K(t) = \sigma(\theta\in S^{n-1}: |K\cap \theta^\perp| \ge t),
$$
and
$$
\Pi_K(t) = \sigma(\theta\in S^{n-1}: |K | \theta^\perp| \ge t).
$$
Here $\theta^\perp = \{x\in \R^n: \langle x,\theta\rangle =0\}$, and $K|\theta^{\bot}$ denotes the orthogonal projection of $K$ onto $\theta^{\bot}$. We write $|A|$ for the $k$-dimensional Lebesgue measure (volume)  of a  set $A \subset \R^n$, where $k$ is the dimension of the minimal flat containing $A$.

It is well known (see e.g. \cite{Ga}) that if $K$ and $L$ are origin-symmetric convex bodies in $\mathbb R^n$ such that
$$|K\cap \theta^\perp|  = |L\cap \theta^\perp|, \quad \mbox{ for all } \theta\in S^{n-1},$$
or
$$|K|\theta^\perp|  = |L| \theta^\perp|, \quad \mbox{ for all } \theta\in S^{n-1},$$
then $K=L$.

Since distribution functions provide much weaker information, one cannot expect similar uniqueness results. But one can ask whether the knowledge of some distribution functions may determine the distribution of the radial function or the support function of the body $K$, or the volume of $K$, or some other information.


In this paper we show that there are origin-symmetric convex bodies $K$ and $L$ in
$\mathbb R^n$, $n\ge 3$, that have the same distribution of the areas of central
sections, but their volumes are different. (In particular, this means
that $K$ and $L$ have different distributions of their radial
functions). In $\R^2$ the latter result is actually true.

We also show that the distribution of the areas of projections of a convex body
in $\mathbb R^n$, $n\ge 2$,
does not determine its volume.

Among positive results, we show that the knowledge of the distribution
function of the areas of  non-central $t$-sections, of the form $K\cap (\theta^\perp + t\theta)$, for every
$t>0$ does determine the distribution of the radial function of $K$.

It is also interesting to look at volume comparison problems involving distribution functions. For example, the celebrated Busemann-Petty problem asks the following question. Let $K$ and $L$ be origin-symmetric convex bodies in $\R^n$ such that
$$|K\cap \theta^\perp|  \le |L\cap \theta^\perp|, \quad \mbox{ for all } \theta\in S^{n-1}.$$
Does it follow that $|K|\le |L|$?
The answer to this question is affirmative if $n\le 4$ and negative if $n\ge 5$; see \cite{Ga} or \cite{K}.
If  we replace the sections by projections, then the corresponding problem is known as the Shephard problem. It has a positive answer in $\R^2$, while in higher dimensions there are counterexamples.

We note, however, that there are certain classes of bodies for which  these two problems do have an affirmative answer in all dimensions.

Here we will consider Busemann-Petty-Shephard type questions for
distribution functions. It is, of course, impossible to obtain an affirmative answer in all dimensions, because of the solutions to the original Busemann-Petty and Shephard problems. But one can ask if there is a positive answer within a smaller class of bodies, or if there are isomorphic versions of the problems.


We show, for example, that if $E$ is a centered ellipsoid and $K$ is
an origin-symmetric convex body such that $S_E(t)\le S_K(t)$ for all
$t$, then $|E|\le |K|$. If $E$ is replaced by a convex intersection
body $L$, this is no longer true. However, in the latter case, there
is an absolute constant $C$ such that $S_L(t)\le S_K(t)$ implies
$|L|\le C |K|$.
We also establish similar results for the case of projections
(intersection bodies would need to be replaced by polar projection bodies).

 The paper is organized as follows: in Section 2 we present all required basic definitions from Convex Geometry and Harmonic Analysis. In Section 3, to develop some intuition, we study two-dimensional cases of our problems. Higher dimensional results are presented in Section 4.
\\

\noindent {\bf Acknowledgment}:  We are indebted to Alexander Koldobsky  and Dmitry Ryabogin for many valuable discussions.

\section{Notation and preliminaries}

Here we will collect some necessary definitions and facts. As usual,   we denote by $\langle x,y\rangle$  the inner product of two vectors $x, y \in \R^n$ and by $|x|$ the length of a vector $x \in \R^n$. We also denote by $B_2^n=\{x \in \R^n: |x| \le 1\}$ the Euclidean unit ball.  Recall that a convex body in $\R^n$ is a convex compact set with non-empty interior. We will assume that the origin is an interior point of $K$.
We say that a compact set $K\subset \R^n$ is a {\it star body} if it is star-shaped about the origin and its {\it radial function} defined by
$$\rho_K(x)=\max \{a\ge 0: ax \in  K \}, \qquad x\in \mathbb R^n\setminus\{0\}, $$ is positive and continuous.

The  {\it Minkowski functional}  of a star body $K$ is
defined by
$$\|x\|_K=\min \{a\ge 0: x \in aK \}, \qquad x\in \mathbb R^n.$$
Of course, $\rho_K(x) = \|x\|_K^{-1}$. We also note that a convex body that contains the origin in its interior is also a star body.

A star body $K$ is origin-symmetric if $K = -K$.
The {\it support function} of a convex body $K\subset \R^n$ is given by
$$h_K(x) = \sup_{y\in K} \langle x,y\rangle, \qquad x\in \mathbb R^n.$$
The {\it polar body} of $K$ is $K^\circ=\{x\in \R^n: \langle x,y\rangle\le 1, \mbox{ for all } y\in K\}.$ Note that $\rho_{K^\circ} = 1/h_K$.

One of the techniques extensively used this paper is the Fourier transform approach developed by Koldobsky; we refer to \cite{K, KY, RZ} for more details, here we will only give basic definitions.
We denote by ${\mathcal S}$ the space of rapidly decreasing infinitely differentiable functions (test functions) on $\R^n$ with values in
${\mathbb C}$. By ${\mathcal S}'$ we denote the space of distributions over
${\mathcal S}$. Every locally integrable real-valued function $f$ on
$\R^n$ with power growth at infinity represents a
distribution acting by integration: for every $\phi \in {\mathcal S},$
$ \langle f, \phi \rangle=\int_{\R^n} f(x) \phi(x)dx.$
The Fourier transform of a distribution $f$ is defined by
$ \langle \hat{f},
\phi\rangle=\langle f,\hat{\phi}\rangle,$
for every test function $\phi$.
 A distribution $f$ is called {\it even
homogeneous} of degree $p \in \R$ if
$$
\langle f(x), \phi(x/t)\rangle=|t|^{n+p}\langle f(x), \phi(x)\rangle,\,\,\,
\forall \phi\in {\mathcal S}, \,\, t\in \R\setminus\{0\}.
$$
The Fourier transform of an even homogeneous distribution of degree $p$ is
an even homogeneous distribution of degree $-n-p$.

A distribution $f$ is called {\it positive definite} if, for every
nonnegative test function  $\phi\in {\mathcal S}$,
$$
\langle f,\phi*\overline{\phi(-x)}\rangle\ge 0.
$$
By Schwartz's generalization of Bochner's theorem, a distribution is positive
 definite if and only if its Fourier transform is a positive distribution
(in the sense that $\langle \hat{f},\phi\rangle\ge 0$, for every non-negative
$\phi\in {\mathcal S}$).

The {\it spherical Radon transform} is a bounded linear operator on
$C(S^{n-1})$ defined by
$$
{\mathcal R} f(\xi)=\int_{S^{n-1}\cap \xi^\bot} f(x) dx,\,\,\, f\in C(S^{n-1}),\,\, \xi\in S^{n-1}.
$$
Koldobsky (\cite{K1}, Lemma 4) proved that if $g(x)$ is an even homogeneous
function of degree $-n+1$ on $\R^n\setminus \{0\}$, $n>1$, so that
$g\big|_{S^{n-1}}\in L_1(S^{n-1})$, then
\begin{equation}\label{eq:kol}
{\mathcal R}g(\xi)=\frac{1}{\pi}\hat{g}(\xi), \,\,\,\, \,\,\,\,\, \forall \xi \in S^{n-1}.
\end{equation}

Let $K$ be an origin-symmetric star body in $\R^n$. Its intersection
body $IK$ is the star body with the radial
function $$\rho_{IK}(\theta) = | K \cap \theta^\perp|, \quad \theta
\in S^{n-1}.$$
Formula  (\ref{eq:kol}) implies that 
\begin{equation}\label{eq:intbodyfourier}
\rho_{IK} =c_1(\rho_{K}^{n-1})^\wedge\quad\text{and}\quad
\rho_{K} =c_2(\widehat{\rho_{IK}})^{\frac1{n-1}},
\end{equation}
where $c_1=[\pi(n-1)]^{-1}$ and $c_2^{-n+1}=(2\pi)^nc_1$.

The concept of intersection bodies of star bodies  was introduced by Lutwak \cite{L} and it played an
important role in the solution of the Busemann-Petty problem,
mentioned in the introduction.

A more general class of bodies is defined as follows. A star body
$K\subset \R^n$ is said to be an intersection body if there exists a
finite Borel measure $\mu$ on $S^{n-1}$ so that $\rho_K =
\RR\mu$. A characterization of intersection bodies via the  Fourier
transform was discovered by Koldobsky as an application of formula (\ref{eq:kol}) (see  \cite[Theorem 4.1]{K}  for more details):

{\it  An origin-symmetric star body $K$ is an intersection body if and only if $\|\cdot\|_K^{-1}$ represents a positive definite distribution on $\R^n$. }

A generalization of the concept of intersection bodies was introduced by Koldobsly; see \cite[Section 4.2]{K}. Let $1\le k \le n-1$, and let $K$ and $L$ be origin-symmetric star bodies in $\R^n$. We say that $K$ is the $k$-intersection body of $L$ if for every $(n-k)$-dimensional subspace $H$ of $\R^n$ we have $$|K\cap H^\perp| = |L \cap H|.$$
Note that not every body $L$ has its $k$-intersection body. But in the case when the $k$-intersection body of $L$ exists, we will denote it by $I_kL$.

Let $K$ be a convex body in $\mathbb R^n$. Define the {\it parallel section function} of $K$ in the direction of $\theta \in S^{n-1}$ by
$$A_{K,\theta }(z)= |K\cap \{\theta^\perp+z\theta\}|,\quad z\in \mathbb R.$$
If $K$ is sufficiently smooth, the fractional derivative of order $q$ of the parallel section function at zero is defined by
$$A_{K,\xi}^{(q)}(0) = \left\langle \frac{t_+^{-1-q}}{\Gamma(-q)}, A_{K,\xi}(t)\right\rangle,$$
where $t_+=\max \{0,t\}$.

The projection body of a convex body $K$ in $\mathbb{R}^n$ is defined
as the body $\Pi K$ with support function
$$h_{\Pi K}(\theta)=|K|\theta^{\bot}|,  \mbox{ for all   }  \theta\in S^{n-1}.$$

There is a Fourier characterization of projection bodies similar to
that for intersection bodies (see \cite{K, KY, RZ}):
{\it An origin symmetric convex body $L$ in $\R^n$ is a projection body  if and only if
$-h_L$ is a positive definite distribution outside of the origin.
}

We finally note that polars of projection bodies are usually called polar projection bodies.

\section{Case of $\R^2$ and first observations}

We will start with some simple observations.
Suppose that we have two origin-symmetric convex bodies $K$ and $L$ in ${\mathbb R}^2$ such that
$$
S_K(t)=S_L(t),  \mbox{  for all   }  t \ge 0.
$$
Since in $\R^2$ the intersection body of $K$ is obtained by rotating $K$ by $\pi/2$ and expanding it by the factor of 2, the condition above is equivalent to
$$
\sigma(\theta:  \rho_{K}(\theta) \ge t)=\sigma(\theta: \rho_{L}(\theta) \ge t),  \mbox{ for all   }  t \ge 0.
$$
Thus, in $\R^2$ the distribution of the areas of central sections does determine the distribution of the radial function. (Later we will see that this is not the case in higher dimensions).

It is easy to see that the distribution of the radial function is not enough to determine the body. Take, for example,
$$
K=B_2^2 \cap \{x: |x_1| \le  1 - \epsilon\} \cap \{x: |x_2| \le 1 - \epsilon\}
$$
and
$$
L=B_2^2 \cap \{x: |x_1| \le 1 - \epsilon\} \cap \{x: |x\cdot (\frac{1}{\sqrt{2}}, \frac{1}{\sqrt{2}})| \le  1 - \epsilon\},
$$
where $\epsilon>0$ is small enough. Then $K$ and $L$ are different, but have the same distribution of the radial functions.

A similar observation allows to conclude that in $\R^3$ (as well as
higher dimensions) convex bodies with the same distribution of radial
functions may have   different distributions of the areas of their central sections.

For $\epsilon>0$ small enough, let
$$
K=B_2^3 \cap \{x: |x_1| \le 1 - \epsilon\} \cap \{x: |x_2| \le 1 - \epsilon\} \cap \{ |x\cdot (\frac{1}{\sqrt{2}}, \frac{1}{\sqrt{2}}, 0)| \le 1 - \epsilon\}
$$
and
$$
L=B_2^3 \cap \{x: |x_1| \le 1 - \epsilon\} \cap \{x: |x_2| \le 1 - \epsilon\} \cap \{ |x\cdot (\frac{1}{\sqrt{3}}, \frac{1}{\sqrt{3}}, \frac{1}{\sqrt{3}})| \le 1 - \epsilon\}.
$$
Then
$$
 \sigma(\theta:\rho_K(\theta) \ge t) = \sigma(\theta:\rho_L(\theta) \ge t),
 $$
but
$$
\sigma(\theta: |K\cap  \theta^\perp| \ge t)\not = \sigma(\theta: |L\cap  \theta^\perp| \ge t).
$$


Let us now move on to volume comparison results.

\begin{proposition}\label{sections_2dim}
Let  $K$ and $L$  be origin-symmetric convex bodies in $  \R^2$ such that
$$
S_K(t)\le S_L(t),  \mbox{  for all   }  t \ge 0.
$$
Then
$$
|K| \le |L|.
$$
\end{proposition}
\bp If
$$
\sigma(\theta: |K\cap \theta^\perp| \ge t)\le \sigma(\theta: |L\cap \theta^\perp| \ge t),  \mbox{  for all   }  t \ge 0,
$$
then we have
$$
\sigma(\theta: 2\rho_K(\theta) \ge t)\le \sigma(\theta: 2\rho_L(\theta) \ge t),  \mbox{  for all    }  t \ge 0,
$$
and so
$$
\sigma(\theta: \rho_K^2(\theta) \ge t)\le \sigma(\theta: \rho_L(\theta)^2 \ge t),  \mbox{  for all   }  t \ge 0.
$$
Integrating the latter inequality over $t \ge 0$ and using the Fubini theorem we get
$$
\int_{S^1} \rho_K^2(\theta) d\theta \le  \int_{S^1} \rho_L^2(\theta) d\theta,
$$
and thus $|K| \le |L|$.
\ep
We note that the  analogous result is not true for projections.
\begin{proposition}\label{projections_2dim}
There exist  origin-symmetric convex bodies $K$, $L \subset \R^2$ such that
$$
\Pi_K(t)= \Pi_L(t),  \mbox{ for all }  t \ge 0,
$$
but
$$
|K|\ne |L|.
$$
\end{proposition}
\bp
Using that in dimension two the length of projections of an origin-symmetric convex body can be written in terms of its support function, we get
$$
\sigma(\theta: h_K (\theta) \ge t)\le \sigma(\theta: h_L(\theta) \ge t),  \mbox{  for all   }  t \ge 0.
$$
It will be convenient to think of support functions as functions on the interval $[0,2\pi]$ and  use the formula    \cite[2.4.27]{Gr}:
$$
|K|  =\frac{1}{2} \int\limits_0^{2\pi} \left(h_K^2(x) -  (h_K'(x))^2\right) dx.
$$
Our goal is to construct two functions $h_K$ and $h_L$  that are the support functions of two symmetric convex bodies, and such that they have equal distribution functions, but
$$
 \int\limits_0^{2\pi} (h_K'(x))^2  dx\not =  \int\limits_0^{2\pi} (h_L'(x))^2 dx.
$$
To do so, we will consider a function $h \in C^2 (\mathbb R)$ that is $\pi$-periodic, and such that  $h(0)=h(\pi)=0$, and $h(x)$ is strictly increasing on $[0, \pi/2]$ and $h(\pi-x)=h(x)$ for $x \in [0, \pi/2]$. Let us denote by $h_1$ the restriction of $h$ to the interval $[0,\pi/2]$ and by $h_2$ its restriction to the interval $[\pi/2,\pi]$. Assume $h(\pi/2)=1$. Then the inverse functions $h_1^{-1}$ and $h_2^{-1}$ are defined on the interval $[0,1]$ with the corresponding ranges $[0,\pi/2]$ and $[\pi/2,\pi]$. Now we have
$$
\int_0^{\pi} (h'(x))^2 dx = \int_0^{\pi/2} (h_1'(x))^2 dx + \int_{\pi/2}^{\pi} (h_2'(x))^2 dx.
$$
Making the change of variables $x=h_1^{-1}(y)$ in the first integral and $x=h_2^{-1}(y)$ in the second, we get
$$
 = \int_0^1 \frac{1}{[h_1^{-1}(y)]'} dy + \int_1^0 \frac{1}{[h_2^{-1}(y)]'} dy.
$$
If we denote $f_1=h_1^{-1}$ and $f_2=\pi -f_1 =h_2^{-1}$, then
$$
 = \int_0^1 \frac{1}{f_1'(y)} dy + \int_1^0 \frac{1}{f_2'(y)} dy=\int_0^1 \frac{2}{f_1'(y)} dy.
$$
Next we will replace the function $h$ by an equally distributed $\pi$-periodic function $\bar{h} \in C^2(\mathbb R)$. For this we take a function $g\in C^2 [0,1]$ such that ${\rm supp}\, g \subset [1/4, 3/4]$ and replace
$f_1$ and $f_2$ with $f_1+\delta g$ and   $f_2+\delta g$, where $\delta >0$ is small enough. Then
$$
\int_0^{\pi} (\bar h'(x))^2 dx = \int_0^1 \frac{1}{f_1'(y)+\delta g'} dy + \int_1^0 \frac{1}{f_1'(y)+\delta g'} dy$$
$$=\int_0^1 \left( \frac{1}{f_1'(y)+\delta g'} + \frac{1}{f_1'(y)-\delta g'} \right)dy
 =\int_0^1 \left( \frac{2 f'_1(y)}{(f_1'(y))^2- (\delta g')^2} \right)dy > \int_0^1 \frac{2}{f_1'(y)} dy.
$$
To finish the proof we set $h_K=1+\eps h$ and $h_L=1+\eps \bar{h}$, where $\eps>0$ is small enough.
\ep

Let us also provide a construction that yields  two infinitely smooth origin-symmetric convex bodies in $\R^2$ whose radial functions have the same distribution function, while the support functions do not. The polars of these  bodies can  be used  to give  another proof of Proposition~\ref{projections_2dim}. In the next section we will also use this construction in higher dimensions. Denote by $e_1$, $e_2$  the standard orthonormal basis of $\R^2$. Let $\eps>0$ be small and consider two bodies $K_0$, $E_0\subset\R^2$ defined by their radial functions, regarded as functions on the interval $[0,2\pi]$,
\begin{equation}\label{eq:radialE0}
\rho_{E_0}(u) = \left(1+\eps\sin^2\bar u\right)^{-1/2},\quad\bar u\in[0,2\pi]
\end{equation}
and
\begin{equation}\label{eq:radialK0}
\rho_{K_0}(u) = \left(1+\eps\sin^22\bar u\right)^{-1/2},\quad\bar u\in[0,2\pi],
\end{equation}
where $\bar u$ is the angle between $u\in S^1$ and    $e_1$.
Then $E_0$ is an  origin-symmetric ellipse in $\R^2$, and $K_0$ is an infinitely smooth origin-symmetric body in $\R^2$ that is not an ellipse. Note that $K_0$ is convex for small enough  $\eps>0$  (see Lemma~\ref{lem:smooth_near_ball}  below). Moreover, the radial functions of $K_0$ and $E_0$ have the same distribution function. Indeed, for each $t>0$,
\begin{align*}
\sigma\left(u\in S^1:\rho_{K_0}(u)\ge t\right) &= \frac1{2\pi}\left|\big\{\bar u\in[0,2\pi]:\rho_{K_0}(u)\ge t\big\}\right| \\
&= \frac2{2\pi}\left|\big\{\bar u\in[0,\pi]:\rho_{K_0}(u)\ge t\big\}\right|.
\end{align*}
Since the set of all $\bar u\in[0,\pi]$ with $\rho_{K_0}(u)\ge t$ is the same as the set of all $\bar u/2\in[0,\pi]$ with $\rho_{E_0}(u)\ge t$, the above measure is equal to
\begin{align*}
\frac2{2\pi}\left|\big\{\bar u/2\in[0,\pi]:\rho_{E_0}(u)\ge t\big\}\right| &= \frac1{2\pi}\left|\big\{\bar u\in[0,2\pi]:\rho_{E_0}(u)\ge t\big\}\right| \\
&=\sigma\left(u\in S^1:\rho_{E_0}(u)\ge t\right).
\end{align*}

Let $K$ and $E$ be the polar bodies of $K_0$ and  $E_0$, respectively.  Since $E$ is an ellipsoid and $K$ is not, the equality case of the Blaschke-Santal\'o inequality implies $$|K||K_0|<|E||E_0|.$$
 We now use that $|K_0| = |E_0|$ to see that $K$ and $E$ have different volumes. Since the radial function of the polar body is equal to the reciprocal of the support function of the original body, the latter implies that the support functions of $K_0$ and  $E_0$ have different distribution functions.

 Moreover, we obtain another proof of Proposition~\ref{projections_2dim}, if we observe that
\begin{align*}
\Pi_K(t) &= \sigma\left(u\in S^1:h_K(u)\ge t/2\right) =\sigma\left(u\in S^1:\rho_{K_0}(u)\le 2/t\right) \\
&= \sigma\left(u\in S^1:\rho_{E_0}(u)\le 2/t\right) =\sigma\left(u\in S^1:h_E(u)\ge t/2\right) = \Pi_E(t).
\end{align*}

\begin{remark} Proposition \ref{projections_2dim} shows that the knowledge of the distribution of the length of projections of a convex body in $\R^2$ does not give the distribution of the radial function of the body.
\end{remark}

Note that if in Proposition \ref{projections_2dim} one of the bodies is a disk, then not only the volumes are the same, but also the second body must be a disk. In fact, this is true in all dimensions.

\begin{proposition} Let  $K$ be an origin-symmetric convex body in $\R^n$ such that for some  $r>0$ we have
$$
\Pi_K(t)=\Pi_{rB_2^n}(t) ,  \mbox{  for all   }  t \ge 0.
$$
Then $K=rB_2^n$.
\end{proposition}
Analogously, for sections:
\begin{proposition} Let  $K$ be an origin-symmetric convex body in $\R^n$ such that for some  $r>0$ we have
$$
S_K(t)=S_{rB_2^n}( t),  \mbox{  for all   }  t \ge 0.
$$
Then $K=rB_2^n$.
\end{proposition}

These facts easily follow  from the observation
$$ \sigma(\theta:|rB_2^n| \theta^\perp| \ge t) =\sigma(\theta:|rB_2^n \cap  \theta^\perp| \ge t) =
\left\{\begin{array}{ll}
   1, & \mbox{    if    } t < r^{n-1}|B_2^{n-1}|,\\
 0, & \mbox{    if    } t \ge  r^{n-1}|B_2^{n-1}|.
  \end{array}\right.
$$

The same observation allows to obtain the following volume comparison result.

\begin{proposition}
Let $L$ be a convex origin-symmetric body in $ \mathbb R^2$ such that
$$
\Pi_{r B_2^2}(t)\le \Pi_L(t), (\mbox{or  } \Pi_{r B_2^2}(t)\ge \Pi_L(t)),  \mbox{ for all }  t \ge 0.
$$
Then $|r B_2^2|\le |L|$, (or $|r B_2^2|\ge |L|$).

\end{proposition}
\bp
The hypothesis of the proposition implies that $|r B_2^2 | \theta^\perp| \le |L | \theta^\perp|$ (correspondingly, $|r B_2^2 | \theta^\perp| \ge |L | \theta^\perp|$) for every $\theta \in S^1$. The positive answer to the Shephard problem finishes the proof.
\ep

For the case of projections we also have the following.
\begin{proposition}
Let  $K,L \subset \R^2$ be origin-symmetric convex bodies such that
$$
\Pi_K(t) \le \Pi_L(t),  \mbox{ for all }  t \ge 0.
$$
Then
$$
|\partial K| \le |\partial L|  \mbox { and } |K^\circ|\ge|L^\circ|.
$$

\end{proposition}

\bp
Note that the hypothesis of the proposition is equivalent to
$$
\sigma(\theta: h_K(\theta) \ge t)\le \sigma(\theta: h_L(\theta) \ge t),  \mbox{ for all }  t \ge 0.
$$
Integrating the latter with respect to $t$, we get $$\int_{S^1} h_K(\theta) d\theta \le \int_{S^1} h_L(\theta) d\theta,$$ which   is equivalent to $|\partial K| \le |\partial L|$.

Further, observe that $$
\sigma(\theta: [h_K(\theta)]^{-2} \le t^{-2})\le \sigma(\theta: [h_L(\theta)]^{-2}  \ge t^{-2}),  \mbox{ for all }  t > 0,
$$
that is
$$
\sigma(\theta: [h_K(\theta)]^{-2} > t )\ge \sigma(\theta: [h_L(\theta)]^{-2}  > t),  \mbox{ for all }  t > 0.
$$
After integration we get
$$\int_{S^1} [h_K(\theta)]^{-2} d\theta \ge \int_{S^1} [h_L(\theta)]^{-2} d\theta,$$
which means $|K^\circ|\ge|L^\circ|$.

\ep

Note that the fact that $\Pi_K(t) \le \Pi_L(t)$,  for all $ t \ge 0$, implies   $|\partial K|\le |\partial L|$, is true in all dimensions.

\section{Case of $\R^n$, $n\ge 3$: Central hyperplane sections and orthogonal projections}

First, we will prove that, unlike in $\R^2$, the distribution of areas of central sections does not determine the volume of the body. We will give two proofs. The second is shorter, but the first proof will be needed later when we discuss derivatives of section functions. We will be using the following lemma.

\begin{lemma}\label{lem:smooth_near_ball}
Let $f:S^{n-1} \to \R$ be an infinitely smooth even function, and $p\ne 0$ a real number. Then there exists $\eps_0>0$ such that for each $\eps\in(0,\eps_0)$ the function $(1+\eps f)^p$ can be considered as each of the following functions:
\begin{enumerate}
\item[(i)] the radial or support function of a convex body in $\R^n$,
\item[(ii)] the radial function of the $k$-intersection body of a convex
  body, for every $1\le k\le n-1$,
\item[(iii)] the support function of the projection body of a convex body.
\end{enumerate}
\end{lemma}

\bp

Let $F(\xi) = (1+\eps f(\xi))^p$, $\xi \in S^{n-1}$,  and denote by $F_q$ its homogeneous extension of degree $q$ to $\mathbb R^n\setminus \{0\}$, i.e. 
$$F_q(x) = |x|^q F\left(\frac{x}{|x|}\right), \quad x \in \mathbb R^n\setminus \{0\}.$$

We claim that $\rho_K = \left[ F_q\right]^{-1/q}$ (assuming $q\ne 0$) is the radial function of a convex body $K$. This is a standard approximation argument.
First of all, observe that $\rho_K$ is close to the radial function of the Euclidean ball in
$C^l(S^{n-1})$ for every $l\in \mathbb N$. Therefore an application of
the well-known formula for the curvature of a planar curve 
\begin{equation}\label{curvature}
\kappa =
\frac{\rho^2+2(\rho')^2-\rho \rho''}{(\rho^2+(\rho')^2)^{3/2}},
\end{equation}
 and
the fact that the Euclidean ball has a strictly positive curvature
imply that $K$ is convex. Part (i) follows.

To prove (ii),  assume that $q$ is an integer, $-n<q<0$. Let us show that  $F_q(x)$ is a positive definite distribution on $\mathbb R^n$ for small $\eps$.
Let $m=\lfloor(n+q)/2\rfloor$ be the largest integer less than or equal to $(n+q)/2$. If $n+q$ is an odd integer, then \cite[Lemma 3.16]{K} gives
\begin{equation*}
\widehat{F_q}(u)=(-1)^m\pi\int_{S^{n-1}\cap u^\perp}\Delta^mF_q(v)\,dv, \quad u\in S^{n-1}. 
\end{equation*}
If $n+q$ is an even integer, then \cite[Lemma 3.1]{Y} gives
\begin{equation*}
\widehat{F_q}(u) =(-1)^m\pi\int_{S^{n-1}}\left(\ln|\langle u,v\rangle| \Delta^m F_q(v) -(n-2)\Delta^{m-1}F_q(v) \right) \,dv.
\end{equation*}

Since $|\cdot |^q$ is a positive definite distribution and since $F_q$ is close to $|\cdot |^q$ in $C^l(S^{n-1})$ for every 
$l\in \mathbb N$, we see that $F_q$ is also a positive definite distribution. 

Let $G(x) = \left[\widehat{F_q}(x)\right]^{1/(n+q)}$. Using formula (\ref{curvature}),  the connection between the Fourier  transform and differentiation, and the fact that $F_q$ is close to $|\cdot |^q$, it's not hard to show that $G$ defines the radial function of a convex body. Thus part (ii) is proved.

The proof of part (iii) of the lemma is similar to that of part (ii); take $q=1$.

\ep

\begin{theorem}\label{different_volume}
There exist  two origin-symmetric convex bodies $K, L \subset \R^n$, $n \ge 3$, such that
$$
S_K(t)= S_L(t),  \mbox{  for all   }  t \ge 0.
$$
but
$$
|K| \ne |L|.
$$
\end{theorem}
\bp
Let $H_2(x)$ be a second degree zonal spherical harmonic. Let $F$ be an even $C^\infty$ function on the sphere that has the same distribution function as $H_2$ and that differs from $H_2$ in some small neighborhood.
Using Lemma~\ref{lem:smooth_near_ball} we can find origin-symmetric convex bodies $K,L\subset\R^n$ whose intersection bodies are defined by their radial functions as follows:
$$\rho_{IK}(x) = 1 +\eps H_2(x) \quad\mbox{and}\quad \rho_{IL} (x) = 1 + \eps F(x), \quad x \in S^{n-1},$$
where $\eps>0$ is small. Then, by construction,
\begin{equation}\label{distr_eq}
\sigma (\theta: \rho_{IL}(\theta) \ge t) = \sigma (\theta: \rho_{IK}(\theta) \ge t) .
\end{equation}
Using \eqref{eq:intbodyfourier}, we can compute the volume of $K$ as follows:
\begin{align*}
|K| &= \frac{1}{n} \int_{S^{n-1}} \rho_K^{n}(x) dx = \frac{c_2}{n} \int_{S^{n-1}}\left( \widehat{|\cdot|^{-1}}(x)+ \epsilon \widehat{H_2}(x)\right)^{\frac{n}{n-1}} dx \\
&=c_3+c_4\eps\int_{S^{n-1}}\widehat{H_2}(x)\,dx+ +c_5\eps^2\int_{S^{n-1}}\big[\widehat{H_2}(x)\big]^2 dx + O(\eps^3),
\end{align*}
where $c_3$, $c_4$, and $c_5$ are positive constants depending on $n$ only. Here the second term should be zero because the Fourier transform of $H_2$ extended to $\mathbb R^n$ with homogeneity $-1$ is $c(n) H_2$ (see \cite{GYY}). Thus
$$|K| = c_3+ c_5\eps^2 \int_{S^{n-1}} \big[\widehat{H_2}(x)\big]^2 dx + O(\eps^3).$$
Similarly,
$$|L| = c_3+ c_5\eps^2 \int_{S^{n-1}} \big[\widehat{F}(x)\big]^2 dx + O(\eps^3).$$
In order to show that $|K|< |L|$, we will prove that $$\int_{S^{n-1}} \left[\widehat{H_2}(x)\right]^2 dx< \int_{S^{n-1}} \left[\widehat{F}(x)\right]^2 dx .$$

Let $\sum_{m=2}^\infty Q_m$ be the spherical harmonic expansion of $F$, where $Q_m$ is a spherical harmonic of degree $m$. Note that the we only have harmonics of even degrees, and the harmonic of order zero is zero, since $$\int_{S^{n-1}} F(x)  dx =\int_{S^{n-1}} H_2(x)  dx =0.$$
By \cite{GYY}, the spherical harmonic expansion of $\widehat{F}$ is $\sum_{m=2}^\infty \lambda_m Q_m$, where $$\lambda_m= \frac{2^{n-1} \pi^{n/2} (-1)^{m/2}\Gamma((m+n-1)/2)}{\Gamma((m+1)/2)}.$$

Then $$\int_{S^{n-1}} \left[\widehat{H_2}(x)\right]^2 dx = \lambda_2^2 \|H_2\|^2_2,$$
and
$$\int_{S^{n-1}} \left[\widehat{F}(x)\right]^2 dx = \sum_{m=2}^\infty \lambda_m^2 \|Q_m\|^2_2.$$

Observe that equality (\ref{distr_eq}) implies that $$\int_{S^{n-1}} \left[{H_2}(x)\right]^2 dx  = \int_{S^{n-1}} \left[{F}(x)\right]^2 dx ,$$
i.e.
$$\|H_2\|^2_2 = \sum_{m=2}^\infty \|Q_m\|^2_2.$$
It remains to prove that
$$\lambda_2^2 \sum_{m=2}^\infty \|Q_m\|^2_2< \sum_{m=2}^\infty \lambda_m^2 \|Q_m\|^2_2.$$

But this follows from the following two facts. First of all,  the $\Gamma$-function is log-convex and so
$$\log\Gamma\left(\frac{n+1}{2}\right) - \log\Gamma\left(\frac{3}{2}\right) < \log\Gamma\left(\frac{m+n-1}{2}\right) - \log\Gamma\left(\frac{m+1}{2}\right),$$
which implies that $|\lambda_2| < |\lambda_m|$ for all $m>2$. And second of all, one can see that there is a non-zero harmonic $Q_k$ of some degree $k>2$ in the expansion of $F$. If for all $m>2$ we had $Q_m=0$, then $F$ would be a spherical harmonic of degree 2. But if $F$ and $H_2$ are quadratic polynomials, then it's impossible for $F$ to be equal to $H_2$ on some open set and not equal to $H_2$ on some other set.

\ep

Here we present another proof of the previous theorem.
\bp
Let $E_0$ and $K_0$ be the infinitely smooth origin-symmetric convex bodies in $\R^2={\rm span}\{e_1,e_2\}$ defined in \eqref{eq:radialE0} and \eqref{eq:radialK0}. Consider the origin-symmetric convex bodies $K$ and $L$ in $\R^n=\R^2\times\R^{n-2}$ whose intersection bodies are equal to
$$IK=K_0\oplus_2 B_2^{n-2} \quad\text{and}\quad IL=E_0\oplus_2 B_2^{n-2},$$ respectively. Here $\oplus_2$ denotes the $\ell_2$-sum of two origin-symmetric convex bodies.  If  $A\subset\R^n$ and $B\subset\R^m$ are origin-symmetric convex bodies, the sum $A\oplus_2B$ is the body in $\R^n\times\R^m$ with the Minkowski functional $\|(x,y)\|_{A\oplus_2B}=\left(\|x\|_A^2+\|y\|_B^2\right)^{1/2}$. To show the existence of   $K$ and $L$ we will use Lemma~\ref{lem:smooth_near_ball}. To this end, let us describe the radial functions of $IK$ and $IL$. We write an arbitrary vector $\theta$ on the sphere $S^{n-1}$ as
\begin{equation*}
\theta=su+\sqrt{1-s^2}v\in S^{n-1},
\end{equation*}
where $u\in\R^2$ and $v\in\R^{n-2}$ are unit vectors, and $0\le s\le1$. Then
\begin{align}\label{eq:radialK1}
\rho_{IK}(\theta) &= \left\|su+\sqrt{1-s^2}v\right\|^{-1}_{K_0\oplus_2B_2^{n-2}} = \Big(s^2\|u\|_{K_0}^2+(1-s^2)|v|^2\Big)^{-\frac12} \notag \\
&= \left(s^2\rho_{K_0}^{-2}(u)+1-s^2\right)^{-1/2}= \left(1+\eps s^2\sin^22\bar u\right)^{-1/2},
\end{align}
where $\bar u$ is the angle between $u\in S^{n-1}\cap\R^2$ and  $e_1$. Similarly,
\begin{equation}\label{eq:radialL1}
\rho_{IL}(\theta)= \left(s^2\rho_{E_0}^{-2}(u)+1-s^2\right)^{-1/2}= \left(1+\eps s^2\sin^2\bar u\right)^{-1/2}.
\end{equation}
It follows from \eqref{eq:radialK1}, \eqref{eq:radialL1}, and Lemma~\ref{lem:smooth_near_ball} that both $K$ and $L$ are  origin-symmetric  convex bodies when $\eps>0$ is small enough. Equations \eqref{eq:radialK1} and \eqref{eq:radialL1} also show that the radial functions of $IK$ and $IL$ have the same distribution. Indeed, for each $t\ge 0$, consider the set $U = \{ \theta\in S^{n-1}:  \rho_{IK}(\theta)\ge t\}.$ Its Haar measure  is equal to
\begin{equation*}
\int_{S^{n-1}} \chi(\theta)\,d\sigma(\theta)=\int_0^1\frac{s(1-s^2)^{\frac{n-4}2}}{|S^{n-1}|}\left(\int_{S^1}\int_{S^{n-3}} \chi\big(su+\sqrt{1-s^2}v\big)\,dv\,du\right)ds,
\end{equation*}
where $S^1=S^{n-1}\cap\R^2$, $S^{n-3}=S^{n-1}\cap\R^{n-2}$, and $\chi$ is the characteristic function of the set $U$. Since the condition  $\chi\big(su+\sqrt{1-s^2}v\big)=1$ is equivalent to $$\rho_{IK}(\theta)= \left(s^2\rho_{K_0}^{-2}(u)+1-s^2\right)^{-1/2}\ge t,$$ we have
\begin{equation*}
\sigma\left(\theta\in S^{n-1}:\rho_{IK}(\theta)\ge t\right) =\int_0^1\sigma\left(u\in S^1:s^2\rho_{K_0}^{-2}(u)+1-s^2\le t^{-2}\right)g(s)\,ds,
\end{equation*}
where $g(s)=|S^1||S^{n-3}||S^{n-1}|^{-1}s(1-s^2)^{(n-4)/2}$. Similarly, replacing $K$, $K_0$ by $L$, $E_0$ gives
\begin{equation*}
\sigma\left(\theta\in S^{n-1}:\rho_{IL}(\theta)\ge t\right) =\int_0^1\sigma\left(u\in S^1:s^2\rho_{E_0}^{-2}(u)+1-s^2\le t^{-2}\right)g(s)\,ds.
\end{equation*}
Since $\rho_{K_0}$ and $\rho_{E_0}$ have the same distribution as shown in the previous section, the above two integrals coincide. Therefore, $\rho_{IK}$ and $\rho_{IL}$ have the same distribution, i.e.   $S_K(t)=S_L(t)$ for each $t\ge 0$.

To compare the volumes of $K$ and $L$, first note that their intersection bodies have the same volume since
\begin{align*}
|IK|&=\frac1n\int_{S^{n-1}}\big[\rho_{IK}(\theta)\big]^n d\theta = \frac1n\int_{S^{n-1}}\big[\rho_{IL}(\theta)\big]^n d\theta=|IL|.
\end{align*}
On the other hand, since $L$ is an ellipsoid, but $K$ is not, the equality case of the Busemann intersection inequality  (see e.g. \cite[Corollary 9.4.5]{Ga})  implies that
$$|IK|\,|K|^{-n+1}< |IL|\,|L|^{-n+1}.$$
We conclude that   $|K|\neq|L|$.

\ep

In a similar way one can deal with projections.

\begin{theorem} \label{different_volume_projections}
There exist  two origin-symmetric convex bodies $K, L \subset \R^n$,   such that
$$
\Pi_K(t)=\Pi_L(t),  \mbox{  for all   }  t \ge 0.
$$
but
$$
|K| \ne |L|.
$$
\end{theorem}
\bp
Let $E_0,K_0\subset\R^2$ be the planar bodies defined in \eqref{eq:radialE0}, \eqref{eq:radialK0}, respectively, and consider the bodies $K_1=K_0\oplus_2 B_2^{n-2}$ and $L_1=E_0\oplus_2 B_2^{n-2}$ in $\R^n=\R^2\times\R^{n-2}$.
As we saw above, the radial functions of $K_1$ and $L_1$ have the same distribution.  Now define $K$ and $L$ to be the origin-symmetric convex bodies in $\mathbb R^n$ whose polar projection bodies are
$K_1$ and $L_1$, respectively. Furthermore, it follows from \eqref{eq:radialK1}  that,  for each $\theta\in S^{n-1}$,
$$
|K|\theta^\perp| =h_{\Pi K}(\theta)=\rho_{K_1}^{-1}(\theta) = \left(s^2\rho_{K_0}^{-2}(u)+1-s^2\right)^{1/2} =\left(1+\eps s^2\sin^22\bar u\right)^{1/2},
$$
where $s\in[0,1]$ is the length of the projection of $\theta$ onto $\R^2={\rm span}\{e_1,e_2\}$, $u\in\R^2$ is the unit direction vector of the projection, and $\bar u$ is the angle between $u$ and $e_1$.  Similarly,
$$
|L|\theta^\perp| =h_{\Pi L}(\theta)=\rho_{L_1}^{-1}(\theta) = \left(s^2\rho_{E_0}^{-2}(u)+1-s^2\right)^{1/2} =\left(1+\eps s^2\sin^2\bar u\right)^{1/2}.
$$
It follows from Lemma~\ref{lem:smooth_near_ball} that such convex bodies $K$ and $L$ exist. Moreover, since $\rho_{K_1}$ and $\rho_{L_1}$ have the same distribution, the volumes of projections of $K$, $L$ have the same distribution function, i.e. $\Pi_K(t)=\Pi_L(t)$ for each $t\ge 0$.

To compare the volumes of $K$ and $L$, first note that the polars of their projection bodies have the same volume since
\begin{align*}
|(\Pi K)^\circ|=|K_1|&=\frac1n\int_{S^{n-1}}\big[\rho_{K_1}(u)\big]^n du\\
&=\frac1n\int_{S^{n-1}}\big[\rho_{L_1}(u)\big]^n du=|L_1|=|(\Pi L)^\circ|.
\end{align*}
Note also that $L$ is an ellipsoid, and $K$ is a non-ellipsoidal convex body when $\eps>0$ is small enough. The equality case of the Petty projection inequality (see e.g. \cite[Theorem 9.2.9]{Ga}) implies that
$$|K|^{n-1}|(\Pi K)^\circ|< |L|^{n-1}|(\Pi L)^\circ|.$$
We conclude that  $|K|\neq|L|$.

\ep

Theorem \ref{different_volume} implies that a version of the
Busemann-Petty problem for distribution functions of the  areas of central sections has a negative answer in $\mathbb R^n$, $n\ge 3$. Below we will prove some results in the positive direction.

\begin{theorem}\label{ellipsoid_sections}
Let $K \subset \R^n$ be a convex body and $E\subset \R^n$ a centered ellipsoid such that
$$
S_E(t) \le S_K(t),  \mbox{  for all   }  t \ge 0.
$$
Then
$$
   |E| \le  |K|.
$$
\end{theorem}
\bp Using the Fubini theorem we get  $|IE| \le |IK|$. On the other hand,   the  Busemann intersection inequality and its equality case (see  \cite[Corollary 9.4.5]{Ga}) imply that
$$
c(n) |E|^{n-1}= |I E| \le |IK| \le c(n) |K|^{n-1}.
$$
\ep

In a similar fashion, using the Petty projection inequality and its equality case, one has the following result.
\begin{theorem}\label{ellipsoid_projections}
Let $K \subset \R^n$ be a convex body and $E\subset \R^n$ a centered ellipsoid such that
$$
\Pi_K(t) \le \Pi_E(t),  \mbox{  for all   }  t \ge 0.
$$
Then
$$
   |K| \le  |E|.
$$
\end{theorem}

  It is well known that any ellipsoid is an intersection body (see \cite{Ga}, \cite{K}). Thus it is  natural to ask if Theorem \ref{ellipsoid_sections}  will still be true when we replace an ellipsoid  with a general intersection body.   Remark \ref{rm:bpint} below shows that this question has a negative answer in $\R^n$, $n\ge 3$.

\begin{remark}\label{rm:bpint}
In the first proof of Theorem \ref{different_volume} both $K$ and $L$ are intersection bodies for $\eps$ small enough.
\end{remark}
\bp
Use Lemma \ref{lem:smooth_near_ball} with $p=1/(n-1)$.
\ep

For general convex bodies, however, an isomorphic version of  Theorem \ref{ellipsoid_sections}  can be proved.  To state the theorem we will need   the notion of  the isotropic constant $L_K$ of a convex  body $K$. We refer to \cite{BGVV, MP} for the definition and properties of $L_K$. We note that  it follows  from F.~John's theorem that    if $K\subset \R^n$ is an origin-symmetric convex body then $L_K \le c\sqrt{n}$. It was proved by  Bourgain \cite{Bo1, Bo2} that $L_n \le c n^{1/4} \log n$ and the $\log n$ factor was later removed by Klartag \cite{Kl}.
\begin{theorem}
Let  $K$ and $L$ be origin-symmetric convex bodies in $\R^n$ such that
$$
S_K(t) \le S_L(t),  \mbox{  for all   }  t \ge 0.
$$
Then there is an absolute constant $c$ such that
$$
   |K| \le c L_K  |L|.
$$
\end{theorem}
\bp It was proved in \cite[Equation 5.3]{MP}  that for any origin-symmetric convex body $K$ we have
$$
c_1\frac{|K|^{\frac{n-1}{n}}}{L_K} \le \left(\int_{S^{n-1}}   |K\cap \theta^\perp|^n d\sigma(\theta) \right)^{1/n} \le c_2 |K|^{\frac{n-1}{n}} 
$$
where $c_1, c_2>0$ are absolute constants.
Thus, the condition of the theorem gives
$$
c_1\frac{|K|^{\frac{n-1}{n}}}{L_K} \le  c_2 |L|^{\frac{n-1}{n}}
$$
or
$$
|K| \le  c_3 \left( L_K \right)^{\frac{n}{n-1}} |L|.
$$
We use $L_K \le c\sqrt{n}$  to finish the proof.

\ep

It is well known that $L_K$ is bounded above by an absolute constant for different classes of convex bodies in $\R^n$, in particular for convex intersection bodies (combine Corollary 3.2 from \cite{MP} and Theorem 3.4 from \cite{KPY}).

\begin{corollary} There exists an absolute constant $C>0$, such that for any convex intersection body $K$ and any origin-symmetric convex body $L$ in $\R^n$ satisfying
$$
S_K(t) \le S_L(t),  \mbox{  for all   }  t \ge 0,
$$
we have
$$
   |K| \le C  |L|.
$$
\end{corollary}

Next we obtain analogous results for projections instead of sections. To state them, we will need  the concept of the volume ratio. The volume ratio of  a  convex body $K\subset\R^n$, denoted by ${\rm vr}(K)$, is defined by $${\rm vr}(K)=\left(\frac{|K|}{|E|}\right)^{1/n}$$ where $E$ is the ellipsoid of maximal volume contained in $K$.

\begin{theorem}
Let  $K$ and $L$ be origin-symmetric convex bodies in $\R^n$ such that
$$\Pi_K(t) \le \Pi_L(t)\quad \text{ for all }  t \ge 0.$$ Then $$|K| \le {\rm vr}(L)|L|.$$
\end{theorem}
\bp
Using the Petty projection inequality and its reverse form (as shown in \cite{A}), we get
\begin{equation}\label{eq:petty_reverse_vr}
\frac{c(n)}{{\rm vr}(K)} \le |K|^{\frac{n-1}n}|(\Pi K)^\circ|^{\frac1n}\le c(n),
\end{equation}
where $c(n)=|B_2^n|^{\frac{n-1}n}|(\Pi B_2^n)^\circ|^{1/n}=|B_2^n|/|B_2^{n-1}|$.

Note that the condition of the theorem gives
\begin{equation*}
|(\Pi K)^\circ|=|B_2^n|\int_{S^{n-1}}|K|u^\perp|^{-n} d\sigma(u) \ge |B_2^n|\int_{S^{n-1}}|L|u^\perp|^{-n} d\sigma(u) =|(\Pi L)^\circ|.
\end{equation*}
Combining the latter with inequalities \eqref{eq:petty_reverse_vr}, we get $$|K|^{\frac{n-1}n}\le {\rm vr}(L)|L|^{\frac{n-1}n}.$$
The  volume ratio of an origin-symmetric convex body cannot exceed $\sqrt{n}$ (see \cite{Ba1}), which completes the proof.

\ep

It is known that the unit balls of finite-dimensional subspaces of $L_1$ have uniformly bounded volume ratios; see \cite[Theorem 2]{BM} and  \cite[9.3]{MS}, or \cite{Ba} for a direct proof. Note that every polar projection body is the unit ball of a finite-dimensional subspace of $L_1$ (see \cite{K} for more details). Thus the volume ratios of polar projection bodies are bounded by an absolute constant, which yields the following result.
\begin{corollary}\label{cor:isoshephard}
There exists an absolute constant $C>0$, such that for any polar projection body $K$ and any origin-symmetric convex body $L$ in $\R^n$ satisfying
$$
\Pi_K(t) \le \Pi_L(t),  \mbox{  for all   }  t \ge 0,
$$
we have
$$
   |K| \le C  |L|.
$$
\end{corollary}
In fact, Corollary~\ref{cor:isoshephard} can be slightly improved by replacing the polar projection body $K$ with a convex intersection body because convex intersection bodies have uniformly bounded volume ratios \cite[Proposition 6.2]{KYZ}.

\section{Case of $\R^n$, $n\ge 3$: Non-central hyperplane sections and derivatives of the parallel section function}

We now look at the distribution functions associated with the parallel section function.

\begin{theorem}\label{t-sections}
Let $K$ and $L$ be  convex bodies in $\mathbb R^n$. Assume that for every $z\in \mathbb R$ and every $t\ge 0$ we have
$$\sigma(\theta: A_{K,\theta}(z) \ge t) = \sigma(\theta: A_{L,\theta}(z) \ge t).$$
Then
$$\sigma(\theta: \rho_K(\theta) \ge t) = \sigma(\theta: \rho_L(\theta) \ge t).$$

\end{theorem}

\bp
By the hypothesis of the theorem we have $$\int_{S^{n-1}} A_{K,\theta}(z) d\theta = \int_{S^{n-1}} A_{L,\theta}(z) d\theta ,$$ for every $z\in \mathbb R$.
For $p >-1$ consider

$$\int_0^\infty z^p \int_{S^{n-1}} A_{K,\theta}(z)\ d\theta\ dz = \frac12 \int_{S^{n-1}} \int_K |\langle x,\theta\rangle |^p \ dx \ d\theta$$
$$ = \frac{1}{2(n+p)} \int_{S^{n-1}} \int_{S^{n-1}} |\langle \theta,\xi\rangle |^p \rho_K^{n+p}(\xi) \ d\xi \ d\theta = c(n,p)  \int_{S^{n-1}}   \rho_K^{n+p}(\xi) \ d\xi . $$
where $c(n,p)$ is a non-zero constant.

Thus for $p>-1$  we have $$\int_{S^{n-1}}   \rho_K^{n+p}(\xi) \ d\xi  = \int_{S^{n-1}}   \rho_L^{n+p}(\xi) \ d\xi. $$
Both sides of the latter equality are analytic functions of $p\in \mathbb C$ and since they coincide for $p>-1$, they must coincide for all other values of $p$.

Suppose both $K$ and $L$ are contained in a ball of radius $R$. Then for all $m\in \mathbb N \cup \{0\}$ we have the equality of the moments:
$$\int_0^R t^m \sigma(\theta: \rho_K(\theta) \ge t)dt = \int_0^R t^m \sigma(\theta: \rho_L(\theta) \ge t) dt.$$

Now the statement follows from the Hausdorff moment problem (which is a consequence of the Weierstrass approximation theorem); cf. \cite[VII.3]{F}.

\ep

One can also consider the distribution function of the derivatives of the parallel section function at zero. In the theorem below all bodies are assumed to be sufficiently smooth, to guarantee the existence of the corresponding derivatives.

\begin{theorem}\label{derivatives}  Let $q\in (-1,n-1)$ and assume that $q$ is not an odd integer. Then the following properties hold.

\begin{enumerate}
\item Let $q\ne n/2-1$. Then there exist origin-symmetric convex bodies $K$ and $L$ in $\R^n$ such that $$\sigma(\theta: A_{K,\theta}^{(q)} (0) \ge t) = \sigma(\theta: A_{L,\theta}^{(q)}(0) \ge t), \quad\mbox{ for all } t\in \R,$$
but $|K|\ne |L|$.

\item Let $q=n/2 -1$. If $K$ and $L$ are origin-symmetric convex bodies in $\R^n$ such that $$\sigma(\theta: A_{K,\theta}^{(q)} (0) \ge t) = \sigma(\theta: A_{L,\theta}^{(q)}(0) \ge t), \quad\mbox{ for all } t\in \R,$$ then $|K| = |L|$.

\end{enumerate}
\end{theorem}
\bp

Note that if $q$ is an odd integer, then the derivative $A_{K,\theta}^{(q)}(0)$ is zero in every direction and thus does not give us any information about the body.
That is why we consider only the cases when $q$ is a not an odd integer.

  The proof of part (1) is similar to the first proof of Theorem \ref{different_volume}, and we will just outline the main steps.
Recall that
\begin{equation}\label{deriv of A}
A_{K,\theta}^{(q)}(0) = \frac{\cos(\pi q/2)}{\pi(n-q-1)}\left(\|\cdot\|_K^{-n+q+1}\right)^\wedge(\theta), \quad \theta \in S^{n-1};
\end{equation}
see \cite[Theorem 3.18]{K}.

Then $$\sigma(\theta: A_{K,\theta}^{(q)} (0) \ge t) = \sigma(\theta: A_{L,\theta}^{(q)}(0) \ge t), \quad\mbox{ for all } t\in \R,$$
is equivalent to
$$\sigma(\theta: \left(\|\cdot\|_K^{-n+q+1}\right)^\wedge(\theta) \ge t) = \sigma(\theta: \left(\|\cdot\|_L^{-n+q+1}\right)^\wedge(\theta) \ge t), \quad\mbox{ for all } t\in \R.$$

Let $H_2(x)$ be a second degree zonal spherical harmonic. Let $F$ be an even $C^\infty$ function on the sphere that has the same distribution function as $H_2$ and that differs from $H_2$ in some small neighborhood.
Define
$$\left(\|\cdot\|_K^{-n+q+1}\right)^\wedge(x) = 1 +\epsilon H_2(x) \quad\mbox{and}\quad \left(\|\cdot\|_L^{-n+q+1}\right)^\wedge  (x) = 1 + \epsilon F(x), \quad x \in S^{n-1},$$
where $\epsilon>0$ is small.

Extending the two previous equalities to $\R^n\setminus\{0\}$ with homogeneity $-q-1$ and inverting the Fourier transforms, we get
$$\rho_K(\theta) = c(n) \left( 1+  c_3(n) \epsilon\widehat{H_2}(\theta)\right)^{1/(n-q-1)},$$
and
$$\rho_L(\theta) = c(n) \left( 1+c_3(n) \epsilon \widehat{F}(\theta)\right)^{1/(n-q-1)},$$
where the constants $c(n)$ and $c_3(n)$ are different from those before.

After we write the volumes of $K$ and $L$, the problem reduces to comparing the integrals  $\int_{S^{n-1}} \left[\widehat{H_2}(x)\right]^2 dx$ and $\int_{S^{n-1}} \left[\widehat{F}(x)\right]^2 dx $. Unlike in Theorem \ref{different_volume}, $H_2$ and $F$ are now homogeneous of degree $-q-1$. By \cite{GYY}, the Fourier transform of a spherical harmonic $H_m$ of degree $m$, extended to $\R^n\setminus\{0\}$ with homogeneity $-q-1$, is $\lambda_m H_m$, where $$\lambda_m = \frac{2^{n-1} \pi^{n/2} (-1)^{m/2}\Gamma((m+n-q-1)/2)}{\Gamma((m+q+1)/2)}.$$
Using again the log-convexity of the $\Gamma$-function, we obtain that $|\lambda_2|< |\lambda_m|$ if $q<n/2-1$, and $|\lambda_2|> |\lambda_m|$ if $q>n/2-1$. In both cases, the conclusion follows as in Theorem \ref{different_volume}.

We now consider case (2) of the theorem, i.e. when $q = n/2-1$ (and $q$ is not odd).
Let $$\alpha = \min_{\theta \in S^{n-1}} (\|\cdot\|_K^{-n/2})^\wedge(\theta) = \min_{\theta \in S^{n-1}} (\|\cdot\|_L^{-n/2})^\wedge(\theta).$$
Then
\begin{equation}\label{subtract_alpha}
\sigma(\theta: (\|\cdot\|_K^{-n/2})^\wedge(\theta) -\alpha \ge t) = \sigma(\theta: (\|\cdot\|_L^{-n/2})^\wedge(\theta) - \alpha \ge t), \, \mbox{ for all } t\ge 0.
\end{equation}
Integrating both sides of (\ref{subtract_alpha}) with respect to $t\ge 0$ and simplifying, we get
$$\int_{S^{n-1}} (\|\cdot\|_K^{-n/2})^\wedge(\theta) d\theta  = \int_{S^{n-1}}(\|\cdot\|_L^{-n/2})^\wedge(\theta) d\theta. $$
On the other hand, multiplying both sides of (\ref{subtract_alpha}) by $t$, integrating over $t\ge 0$, and using the previous equality, we get
$$\int_{S^{n-1}}\left[ (\|\cdot\|_K^{-n/2})^\wedge(\theta)\right]^2 d\theta  = \int_{S^{n-1}} \left[ (\|\cdot\|_L^{-n/2})^\wedge(\theta)\right]^2 d\theta. $$
Now the spherical Parseval formula (see \cite[Lemma 3.22]{K}) yields
$$\int_{S^{n-1}} \|\theta\|_K^{-n}  d\theta  = \int_{S^{n-1}}  \|\theta\|_L^{-n}  d\theta, $$
that is $|K|=|L|$.

\ep

The following is a consequence of the previous theorem. Let $1\le k\le n-1$ and let $K$ and $L$ be origin-symmetric convex bodies in $\mathbb R^n$ such that their $k$-intersection bodies $I_k K$ and $I_k L$ exist.
Assume $$\sigma(\theta: \rho_{I_k K} \ge t) = \sigma(\theta: \rho_{I_k L}\ge t), \mbox{ for all } t\ge 0.$$
Since $\|\theta\|_{I_k K}^{-k} = c(n,k) \left( \|\cdot \|_K^{-n+k}\right)^\wedge(\theta)$, the previous theorem shows that the volumes of $K$ and $L$ are not necessarily equal   if $k\ne n/2$. However, if $k = n/2$, then  $|K|=|L|$. Moreover, if $k = n/2$, the same ideas can be used to show that $$\sigma(\theta: \rho_{I_k K} \ge t) \le \sigma(\theta: \rho_{I_k L}\ge t), \quad t\ge 0,$$ implies $|K|\le |L|$.

In contrast to Theorem \ref{derivatives}, we have the following.

\begin{theorem}
Let $K$ and $L$ be origin-symmetric convex bodies in $\R^n$ such that
\begin{equation}
\label{deriv}\sigma(\theta: A_{K,\theta}^{(q)} (0) \ge t) = \sigma(\theta: A_{L,\theta}^{(q)}(0) \ge t),
\end{equation}
for all $t\in \mathbb R$ and all $q$ from some interval. Then $$\sigma(\theta: \rho_K(\theta) \ge t)=\sigma(\theta: \rho_L(\theta) \ge t), \mbox{ for all } t\ge 0.$$
\end{theorem}
\bp
Integrating (\ref{deriv}) and using formula (\ref{deriv of A}), we get
$$\int_{S^{n-1}}   \rho_K^{n-q-1}(\xi) \ d\xi  = \int_{S^{n-1}}   \rho_L^{n-q-1}(\xi) \ d\xi.$$
Since the latter integrals are analytic functions of $q\in \mathbb C$ and they coincide on some interval, the latter equality holds for all $q\in \R$. We finish as in Theorem \ref{t-sections} by using the  Hausdorff moment problem.

\ep

\section{Case of $\R^n$, $n\ge 3$: Sections by subspaces of higher co-dimension.}

In conclusion,  we note that one can also consider the distribution of central sections of dimension $k$ for origin-symmetric convex bodies in $\R^n$. Let $\sigma$ be the Haar probability measure on the Grassmanian $Gr(n,k)$ of $k$-dimensional subspaces of $\R^n$. We have the following generalizations of  Theorem~\ref{different_volume} and Theorem~\ref{different_volume_projections} for sections of dimension $2\le k\le n-1$.

\begin{theorem}
Let $2\le k\le n-1$. The following properties hold.
\begin{enumerate}
\item[(i)] Let $K$ be a convex body  and $E$ a centered ellipsoid in $\mathbb R^n$ such that
$$\sigma(H\in Gr(n,k) : |E\cap H| \ge t) \le \sigma(H\in Gr(n,k): |K\cap H| \ge t), \quad \forall t\ge 0.$$
Then $$|E|\le |K|.$$
\item[(ii)]  There are origin-symmetric convex bodies $K$ and $L$  in $\mathbb R^n$ such that
$$\sigma(H\in Gr(n,k): |K\cap H| \ge t) = \sigma(H\in Gr(n,k): |L\cap H| \ge t), \quad \forall t\ge 0,$$
but $$|K|\ne |L|.$$

\end{enumerate}
\end{theorem}

\bp
Part (i) follows from a more general version of the Busemann intersection inequality; see  \cite[p. 372]{Ga}.

To prove part (ii), it suffices to find an origin-symmetric ellipsoid
$L$ and a non-ellipsoidal origin-symmetric convex body $K$ in $\R^n$
that have the same distribution of the areas of  their $k$-dimensional sections, and then use the equality case of the general version of the Busemann intersection inequality.

Let $E_0$ and $K_0$ be the planar bodies in $\R^2={\rm span}\{e_1,e_2\}$ defined as in \eqref{eq:radialE0} and \eqref{eq:radialK0}. Let  $\R^{k+1}$ be a subspace of $\R^n$ containing the fixed subspace $\R^2$. As shown in the second proof of Theorem~\ref{different_volume}, for $\eps>0$ small enough we can find origin-symmetric convex bodies $\bar K$ and $\bar L$ in $\R^{k+1}=\R^2\times\R^{k-1}$ whose intersection bodies are equal to
$$
I{\bar K}=K_0\oplus_2 B_2^{k-1} \quad\text{and}\quad I{\bar L}=E_0\oplus_2 B_2^{k-1}.
$$
As shown in \eqref{eq:radialK1} and \eqref{eq:radialL1}, their radial functions are given by
\begin{align*}
\rho_{I\bar K}(\theta) &= \left(s^2\rho_{K_0}^{-2}(u)+1-s^2\right)^{-1/2}= \left(1+\eps s^2\sin^22\bar u\right)^{-1/2},\\
\rho_{I\bar L}(\theta) &= \left(s^2\rho_{E_0}^{-2}(u)+1-s^2\right)^{-1/2}= \left(1+\eps s^2\sin^2\bar u\right)^{-1/2}.
\end{align*}
where $s\in[0,1]$ is the length of the projection of $\theta$ onto $\R^2$, $u\in\R^2$ is the unit direction vector of the projection, and $\bar u$ is the angle between $u$ and $e_1$.

Note that $I\bar K$ is rotationally invariant with respect to the second coordinate in the sense that $\|x+y\|_{I\bar K}=\|x+\tilde y\|_{I\bar K}$ for each $x\in\R^2$ and $y,\tilde y\in\R^{k-1}$ with $|y|=|\tilde y|$. Note also that the Fourier transform preserves the property of rotational invariance. Thus \eqref{eq:intbodyfourier} implies that $\bar K$ is also rotationally invariant in the second coordinate. So is $\bar L$, due to the same argument.

Now we define two origin-symmetric convex bodies $K$ and $L$ in $\R^n=\R^2\times\R^{n-2}$ by $$\|x+y\|_K=\|x+\tilde y\|_{\bar K}$$ and $$\|x+y\|_L=\|x+\tilde y\|_{\bar L}$$ for each $x\in\R^2$, $y\in\R^{n-2}$, and $\tilde y\in\R^{k-1}$ with $|\tilde y|=|y|$. They are well defined due to the rotational invariance with respect to the second coordinate of $\bar K$ and $\bar L$.

To compute the distribution functions with respect to the Haar measure $\sigma$ on $Gr(n,k)$, fix a subspace $H$ of $\R^n$ of dimension $k$. We may assume that the intersection of  $H$ and $\R^2$ is always one-dimensional because the set of all $H\in Gr(n,k)$ with $\dim (H\cap\R^2)\neq1$ is $\sigma$-null. Indeed, if the dimension of $H\cap\R^2$ is not equal to one, then the intersection is zero- or two-dimensional; each of the corresponding sets is $\sigma$-null as follows:
\begin{equation*}
\sigma(H:\dim (H\cap\R^2)=0)=\sigma(H\in Gr(n,k):H\le\R^{n-2})=0,
\end{equation*}
 and
\begin{align*}
\sigma(H:\dim (H\cap\R^2)=2)&=\sigma(H\in Gr(n,k):H\ge\R^2)\\
&=\sigma(H\in Gr(n,n-k):H\le\R^{n-2})=0.
\end{align*}

Since $H\cap\R^{n-2}$ is $(k-1)$-dimensional, we can choose $\phi\in O(n)$ that maps $H\cap\R^{n-2}$ to $\R^{k-1}$ and fixes each element of $\R^2$. Then $\phi H$ is a subspace of $\R^{k+1}=\R^2\times\R^{k-1}$ of dimension $k$, and it satisfies
\begin{equation*}
\phi(K\cap H)=\bar K\cap\phi H,
\end{equation*}
because $\phi(x+y)=x+\phi y$ and $\|x+y\|_K=\|x+\phi y\|_{\bar K}$ for each $x\in\R^2$, $y\in\R^{n-2}$ with $x+y\in H$.

Let $\theta_H$ be a unit vector in $\R^{k+1}$ perpendicular to $\phi H$. Then the projection of $\theta_H$ onto $\R^2$ can be written as $s_Hu_H$ where $u_H$ is a unit vector in $\R^2$ perpendicular to the line $H\cap\R^2$, and $s_H$ is the distance of $H$ to $\R^2$ in the sense of $$s_H=d(H,\R^2)=\inf\left\{\sqrt{1-\langle\theta,e_1\rangle^2-\langle\theta,e_2\rangle^2}:\theta\in H\cap S^{n-1}\right\}.$$
Thus
\begin{equation}\label{eq:k_sections_dist_K}
|K\cap H| = |\bar K\cap\phi H|=\rho_{I\bar K}(\theta_H)=\left(s_H^2\rho_{K_0}^{-2}(u_H)+1-s_H^2\right)^{-1/2}
\end{equation}
and, similarly,
\begin{equation}\label{eq:k_sections_dist_L}
|L\cap H| = |\bar L\cap\phi H|=\rho_{I\bar L}(\theta_H)=\left(s_H^2\rho_{E_0}^{-2}(u_H)+1-s_H^2\right)^{-1/2}.
\end{equation}

Now we will prove that the distribution function of the volumes of $k$-dimensional sections of $K$ is equal to that of the radial function of $K_0$. To this end, note that the volume of the section of  $K$ by $H$ is equal to the radial function of $I\bar K$ at $\theta_H$ due to \eqref{eq:k_sections_dist_K}, and the interval $(t,1]$ can be rewritten as $$(t,1]=\left\{\rho_{K_0}(u):u\in S^1 \text{ with } \rho_{K_0}(u)>t\right\}.$$
Thus the condition of $|K\cap H|>t$ is equivalent to $$\rho_{I\bar K}(\theta_H)\in\left\{\rho_{K_0}(u):u\in S^1 \text{ with } \rho_{K_0}(u)>t\right\}.$$
 For each measurable set $A\subset S^1$, define the measure $\mu$ on $S^1$ by
\begin{equation*}
\mu(A)=\sigma\Big(H\in Gr(n,k): \rho_{I\bar K}(\theta_H)\in\{\rho_{K_0}(u):u\in A\} \Big).
\end{equation*}
It is not hard to see that it is a regular Borel measure on $S^1$. We leave the details to the reader. Moreover it is a rotationally invariant probability measure. First of all, $\mu(S^1)=1$, since
$$\rho_{I\bar K}(\theta_H)=\left(s_H^2\rho_{K_0}^{-2}(u_H)+1-s_H^2\right)^{-\frac12}\ge \rho_{K_0}(u_H) \ge \min_{u\in S^1}\rho_{K_0}(u).$$
Also, for any rotation $\psi$ on $\R^2$,
\begin{align*}
\mu(\psi A)&=\sigma\Big(H\in Gr(n,k): \rho_{I\bar K}(\theta_H)\in\{\rho_{K_0}(u):u\in \psi A\} \Big) \\
&=\sigma\Big(\tilde\psi H\in Gr(n,k): \rho_{I\bar K}(\theta_H)\in\{\rho_{K_0}(u):u\in A\} \Big) \\
&=\sigma\Big(H\in Gr(n,k): \rho_{I\bar K}(\theta_H)\in\{\rho_{K_0}(u):u\in A\} \Big) =\mu(A),
\end{align*}
where $\tilde\psi\in O(n)$ is an extension of $\psi$ fixing $\R^{n-2}$. Thus $\mu$ is equal to the unique Haar measure $\sigma$ on $S^1$, i.e.,
\begin{equation*}
\sigma\big(H\in Gr(n,k): |K\cap H|>t\big) = \sigma\big(u\in S^1: \rho_{K_0}(u)>t\big).
\end{equation*}
Similarly, replacing $K$, $K_0$ by $L$, $E_0$ gives
\begin{equation*}
\sigma\big(H\in Gr(n,k): |L\cap H|>t\big) = \sigma\big(u\in S^1: \rho_{E_0}(u)>t\big).
\end{equation*}
Lastly, the fact that $\rho_{K_0}$ and $\rho_{E_0}$ have the same distribution function completes the proof.

\ep

\end{document}